\documentclass[11pt]{article}

\usepackage[top=50mm, bottom=50mm, left=50mm, right=50mm]{geometry}
\usepackage{lineno}
\usepackage{amssymb}
\usepackage{amsmath}
\usepackage{amsthm}
\usepackage{dsfont}
\usepackage{epsfig}
\usepackage{graphicx}
\usepackage{graphics}
\usepackage{float}
\usepackage{subfigure}
\usepackage{multirow}
\usepackage{color}
\usepackage{lineno}
\usepackage{fullpage}
\usepackage[normalem]{ulem} 
\usepackage{makeidx}
\usepackage{xspace}
\usepackage{wrapfig}
\usepackage{mathtools}
\usepackage{authblk}
\usepackage{etoolbox}
\usepackage{truncate}
\usepackage{xkeyval}
\usepackage{xstring}
\usepackage{xcolor}
\usepackage{soul}
\sethlcolor{yellow}

\usepackage{hyperref}
\hypersetup{
colorlinks = true,
linkcolor = {red},
urlcolor = {blue},
citecolor = {green}
}

\usepackage[final]{changes}
\definechangesauthor[name={Quan}, color=orange]{quan}

\makeindex

\newtheorem{theorem}{Theorem}
\newtheorem{conjecture}[theorem]{Conjecture}
\newtheorem{definition}[theorem]{Definition}
\newtheorem{corollary}[theorem]{Corollary}
\newtheorem{proposition}[theorem]{Proposition}
\newtheorem{remark}[theorem]{Remark}
\newtheorem{lemma}[theorem]{Lemma}

\newtheorem{observation}[theorem]{Observation}

\newcommand{\Ebb}{\mathbb{E}}

\newcommand{\Rbb}{\mathbb{R}}

\newcommand{\Nbb}{\mathbb{N}}

\newcommand{\Zbb}{\mathbb{Z}}
\newcommand{\Pbb}{\mathbb{P}}

\newcommand{\Hcal}{\mathcal{H}}
\newcommand{\Scal}{\mathcal{S}}

\newcommand{\Lcal}{\mathcal{L}}
\newcommand{\Ocal}[1]{\mathcal{O}\left(#1\right)}
\newcommand{\Xcal}{\mathcal{X}}

\newcommand{\drm}{\textrm{d}}

\newcommand\iid{i.i.d. }

\newcommand{\ie}{i.e. }
\newcommand{\eg}{e.g. }

\newcommand{\ddt}{\dfrac{d}{dt}}
\newcommand{\esperance}[1]{\Ebb\left[ #1 \right]}
\newcommand{\variance}[1]{\textrm{Var}\left[ #1 \right]}

\newcommand{\variancewithstartingpointmodified}[2]{\overline{\textrm{Var}}_{#1}\left[ #2 \right]}

\newcommand{\esperancewithstartingpoint}[2]{\Ebb_{#1}\left[ #2 \right]}
\newcommand{\esperancewithstartingpointmodified}[2]{\overline{\Ebb}_{#1}\left[ #2 \right]}
\newcommand{\norm}[1]{\left\| #1 \right\|}

\newcommand{\proba}[1]{\Pbb\left[ #1 \right]}
\newcommand{\probawithstartingpoint}[2]{\Pbb_{#1}\left[ #2 \right]}
\newcommand{\probawithstartingpointmodified}[2]{\overline{\Pbb}_{#1}\left[ #2 \right]}
\newcommand{\indicator}[1]{\mathds{1}_{\left\{#1\right\}}}

\newcommand{\tmixetaepsilon}{\textrm{t}_\textsc{mix}(\eta;\epsilon)}

\newcommand{\crochet}[2]{\left<#1, #2\right>}
\newcommand{\restr}[2]{{% we make the whole thing an ordinary symbol
  \left.\kern-\nulldelimiterspace % automatically resize the bar with \right
  #1 % the function
  \vphantom{\big|} % pretend it's a little taller at normal size
  \right|_{#2} % this is the delimiter
  }}

\renewcommand{\phi}{\varphi}

\renewcommand{\varpi}[1]{\textrm{Var}_\pi\left[#1\right]}
\newcommand{\tmix}{\textrm{t}_\textsc{mix}(\epsilon)}
\newcommand{\minset}[1]{\inf \left\{ #1  \right\}}
\newcommand{\maxset}[1]{\max \left\{ #1  \right\}}

\newcommand{\dtv}[2]{\textrm{d}_\textsc{tv}\left( #1, #2 \right)}

\newcommand{\wrt}{w.r.t }

\newcommand{\absolutevalue}[1]{\left|#1\right|}

\newcommand{\orangeauthor}[1]{\textcolor{black}{#1}}
\newcommand{\timeincrossinglemma}{t}
\newcommand{\ber}[1]{\textrm{Ber}(#1)}
\newcommand{\mutilde}{\tilde{\mu}}
\newcommand{\densityleft}{p}
\newcommand{\densityright}{q}

\newcommand{\exponentialconstant}{a}
\newcommand{\symmetricgroup}{\Scal_N}
\newcommand{\greentrajectory}{\Hcal_G}

\newcommand{\setcolorinterchange}{\Xcal}
\newcommand{\redbluettwo}{R(X_{t_2}) \cup B(X_{t_2})}

\newcommand{\xisamplingone}{\Xi^1}
\newcommand{\xisamplingn}{\Xi^N}
\newcommand{\pushforwardfunction}{f_*}
\newcommand{\zetatilde}{\tilde{\zeta}}
\newcommand{\zetaredblue}{\zeta^{RB}_{t_2}}
\newcommand{\zetaredbluetilde}{\zetatilde^{RB}_{t_2}}
\newcommand{\identity}{\textrm{Id}}
\newcommand{\redvector}{(R, \dots, R)}
\newcommand{\sigmatilde}{\tilde{\sigma}}
\begin{document}

\title{\LARGE Cutoff for the non reversible SSEP with reservoirs}
\author{\Large Hong Quan Tran\thanks{tran@ceremade.dauphine.fr}}
\affil{\large CEREMADE, Université Paris-Dauphine, PSL University}
\maketitle
%\linenumbers
% \setcounter{page}{0}

\begin{abstract}
    We consider the Symmetric Simple Exclusion Process (SSEP) on the segment with two reservoirs of densities $p, q \in (0,1)$ at the two endpoints. We show that the system exhibits cutoff with a diffusive window, thus confirming a conjecture of Gantert, Nestoridi, and Schmid in \cite{Gantert2020}. In particular, our result covers the regime $p \neq q$, where the process is not reversible and there is no known explicit formula for the invariant measure. Our proof exploits the information percolation framework introduced by Lubetzky and Sly, the negative dependence of the system, and an anticoncentration inequality at the conditional level. We believe this approach is applicable \orangeauthor{to other }models.
\end{abstract}
\tableofcontents

\section{Introduction}
\subsection{Model}
The simple exclusion process is an interacting particle system where the particles attempt to perform simple random walk on a graph, but they are not allowed to jump on top of each other (the exclusion rule). Since its introduction by Spitzer as a simplified model for a gas of interacting particles in \cite{Spitzer1970} (see also \cite{Liggett2005}), it has been shown to exhibit many interesting phenomena, thus received a lot of attention from mathematicians and theoretical physicists, see, for example, \cite{Kipnis1999, Liggett1999, Wood2020}. Recently, a huge amount of work has been devoted to studying the convergence to equilibrium of the conservative (without reservoirs) model in finite volume, see \cite{Wilson2004, Morris2006a, Lacoin2016, Lacoin2016a, Lacoin2017, Lee1998, Yau1997, Labbe2019, Labbe2020}. In this article, we study the non-conservative variant: the SSEP on the segment (the bulk) with reservoirs at the two endpoints, where the particles are allowed to enter or exit the bulk through the reservoirs. We refer the readers to the papers \cite{Bertini2003, Landim2008} for an introduction and motivations on the model and to \cite{Gantert2020, Goncalves2021, Salez2022} for recent developments. More precisely, let $N \in \Zbb_+$ be the length of the segment, and let $p, q \in [0,1]$ be the densities of the reservoirs at the two endpoints. We consider the process $(\eta_t)_{t\geq 0}$ taking values in the state space $\Omega = \{0,1\}^{N}$, whose infinitesimal generator $L$ acts on an observable $\phi: \Omega \to \Rbb$ by 
\begin{equation}\label{generator_of_SSEP}
\begin{split}
    L\phi(\eta) = &\sum_{i = 1}^{N-1} N^2(\phi(\eta^{i \leftrightarrow i+1}) - \phi(\eta)) \\
    &+ N^2 [\densityleft\phi(\eta^{1,1}) + (1-\densityleft)\phi(\eta^{1,0}) - \phi(\eta)]\\
    &+ N^2 [\densityright\phi(\eta^{N,1}) + (1-\densityright)\phi(\eta^{N,0}) - \phi(\eta)],
\end{split}
\end{equation}
where $\eta^{i\leftrightarrow i+1}, \eta^{i, 0}, \eta^{i, 1}$ are the configurations obtained from $\eta$ by swapping the $i$-th and $(i+1)$-th coordinates, resetting the $i$-th coordinate to $0$, resetting the $i$-th coordinate to $1$, respectively. \orangeauthor{Here, time is accelerated by a factor $N^2$ so that the process is observed on a diffusive time scale.} Then $(\eta_t)_{t\geq 0}$ is a SSEP on the segment $[N]:= \{1,\dots, N\}$ with one reservoir of density $p$ placed at site $1$ and the other reservoir of density $q$ placed at site $N$. If $(p,q) \notin \{(0,0), (1,1)\}$, then the process is irreducible, and there is a unique invariant distribution $\pi$. When $p \neq q$, there is no known explicit formula for $\pi$, and the process is not reversible. The classical theory of Markov processes ensures that the system will converge to the invariant distribution $\pi$, no matter from which initial configuration it starts. The distance to equilibrium is measured with respect to the total variation distance $\dtv{\cdot}{\cdot}$, defined by
\[\dtv{\mu}{\nu} = \max\limits_{A\subset \Omega}|\mu(A) - \nu(A)|, \] 
for any distributions $\mu, \nu$ on $\Omega$. The worst-case distance to equilibrium at time $t$ is defined by 
\begin{equation*}
    d(t) = \max_{\eta \in \Omega} \dtv{\Pbb_\eta\left[\eta(t)\in \cdot \right]}{\pi},
\end{equation*}
where $\Pbb_\eta$ is the law of the process starting from $\eta$.
The speed of convergence is quantified by the so-called mixing times: 
\begin{equation*}
    \tmix = \minset{t \geq 0: d(t) \leq \epsilon}.
\end{equation*} 
\paragraph{Cutoff phenomenon.} Consider a family $(P_n)_{n\geq 1}$ of Markov processes. To lighten the notation, we keep the dependence on $n$ implicit as much as possible. The family $(P_n)_{n\geq 1}$ is said to exhibit \textit{cutoff} if, \orangeauthor{in the limit where $n$ tends to infinity}, the asymptotic behavior of $\tmix$ does not depend on $\epsilon$ anymore: 
\begin{equation*}
    \forall \epsilon \in (0, 1)\, \text{fixed}, \;\dfrac{\tmix}{\textrm{t}_\textsc{mix}} \to 1,
\end{equation*}
with $\textrm{t}_\textsc{mix} = \textrm{t}_\textsc{mix}(1/4)$. This means that $d(t)$ undergoes a phase transition around $\textrm{t}_\textsc{mix}$, where it drops from near $1$ to near $0$ in a time of order $o(\textrm{t}_\textsc{mix})$. When cutoff occurs, a natural question is to determine the window in which the phase transition occurs. \deleted{We speak of the \textit{cutoff window}. }More precisely, the family is said to exhibit cutoff with a \textit{window} of size $\Ocal{\omega_n}$ if $\omega_n = o(\textrm{t}_\textsc{mix})$ and 
\begin{align*}
    \lim_{\alpha \to -\infty} \liminf_{n \to \infty} d(\textrm{t}_\textsc{mix} + \alpha\omega_n) &= 1,\\
    \lim_{\alpha \to \infty} \limsup_{n \to \infty} d(\textrm{t}_\textsc{mix} + \alpha\omega_n) &= 0.
\end{align*}
\orangeauthor{The} cutoff phenomenon \orangeauthor{was} discovered by Aldous, Diaconis, and Shahshahani when studying card shuffling \cite{Diaconis1981,Aldous1986, Diaconis1996}, see also \cite{Levin2017} for an introduction to the subject.

\paragraph{Previous works.} The conservative \orangeauthor{SSEP (\textit{without} reservoirs)} has been thoroughly studied in \cite{Morris2006a, Wilson2004, Lacoin2016, Lacoin2016a, Lacoin2017}. In particular, cutoff, and even the limit profile, have been proved by Lacoin in \cite{Lacoin2016, Lacoin2016a, Lacoin2017} for the segment and the circle. \orangeauthor{On the contrary}, only a few works have been \orangeauthor{written} on the non-conservative model. We mention here some \orangeauthor{recent} development. In \cite{Gantert2020}, Gantert, Nestoridi, and Schmid prove a pre-cutoff for the model: for any fixed $\epsilon \in ]0,1[$,
\begin{equation*}
    \dfrac{1}{2\pi^2} \leq \liminf_{N\to \infty}\dfrac{\tmix}{\log N} \orangeauthor{\leq \limsup_{N\to \infty}\dfrac{\tmix}{\log N}} \leq C, 
\end{equation*}
for some constant $C$ independent of $\epsilon$, and they conjecture that the system exhibits cutoff \orangeauthor{with the right estimate on the mixing time given by the lower bound}. Their proof relies on an extension of the coupling used by Lacoin in \cite{Lacoin2016} for the conservative model. In \cite{Goncalves2021}, Gonçalves, Jara, Marinho, and Menezes study \orangeauthor{the reversible case where the} two reservoirs have the same density, $p = q$. Using Yau's famous relative entropy method in \cite{Yau1991}, they prove that
\begin{equation*}
    \tmix = \dfrac{\log N + \mathcal{O}_\epsilon(1)}{2\pi^2}.
\end{equation*}
Recently, Salez studies in \cite{Salez2022} the model on general graphs where reservoirs can be placed at \orangeauthor{arbitrary} sites. By exploiting the negative dependence property of the system, he proves that, under some mild conditions on the graph,
\begin{equation*}
    \tmix = \dfrac{\log N + \mathcal{O}_\epsilon(1)}{\lambda_N},    
\end{equation*}
where $N$ is the size of the graph and $\lambda_N$ is the spectral gap of the random walk of a single particle. We stress that the works above provide a more comprehensive study than the results we just mentioned. \cite{Gantert2020} is more devoted to studying the case where the random \orangeauthor{motion} of the particles is asymmetric. \cite{Goncalves2021} provides the convergence profile for the system from any smooth initial condition. \cite{Salez2022} is more concerned with the characterization of cutoff. However, as far as we know, cutoff has been established in the works \cite{Salez2022, Goncalves2021} only for the case where every reservoir has the same density, which subsequently implies that the invariant distribution is a product measure, and \orangeauthor{that} the system is reversible. 
\paragraph{Our contribution.}
In this paper, we prove cutoff for the SSEP on the segment with reservoirs of arbitrary densities $p, q \in (0, 1)$, thus confirming the conjecture of Gantert, Nestoridi, and Schmid (see Conjecture 1.7 in \cite{Gantert2020}), and also making a \deleted{first }step \orangeauthor{towards} the study of irreversible models. Our proof exploits the information percolation framework introduced by Lubetzky and Sly in \cite{Lubetzky2016}, the negative dependence of the system, and an anticoncentration inequality at the conditional level. We believe that this approach is applicable \orangeauthor{to other }models. 
\subsection{Results}
Our main result is that \orangeauthor{the} model exhibits cutoff at time $\dfrac{\log N}{2\pi^2}$ \orangeauthor{with a window of order $\Ocal{1}$}, as conjectured by Gantert, Nestoridi, and Schmid. \deleted{Moreover, the cutoff window is of order $\Ocal{1}$.}
\begin{theorem}[Main theorem]\label{maintheorem}
    For any $p, q, \epsilon \in (0,1)$ fixed, \orangeauthor{we have}\deleted{ the system exhibits cutoff:}
    \begin{equation*}
        \tmix = \dfrac{\log N}{2\pi^2} + \mathcal{O}_{p,q,\epsilon}(1).
    \end{equation*}
\end{theorem}
In fact, we will prove a stronger \orangeauthor{result}, which makes precise the dependence of the lower order term on $p, q, \epsilon$ and is subsequently still valid when we allow $p, q$ to vary with $N$. Without loss of generality (by the symmetry between $p$ and $q$ and the duality between particles and holes), we suppose that
\begin{equation}\label{assumption}
    q \leq \min\{p, 1-p\}.
\end{equation}
We \orangeauthor{define the weight of the configuration $\eta$ to be}
\[S(\eta) = \sum_{i\in [N]} \eta(i).\]
We denote by $\Ebb^0$ the expectation \wrt the \orangeauthor{absorbing} model where \orangeauthor{the} two reservoir densities $p, q$ are zero. Let $\mathds{1}$ denote the configuration where every site is occupied, and let $t^*$ be the time that the expected weight of the process starting from $\mathds{1}$ falls under a specific threshold:  
\begin{equation}\label{t_star_definition}
    t^* := \minset{t \geq 0: \Ebb_{\mathds{1}}^0\left[S(\eta_t)\right] \leq \sqrt{Np} \vee 1}.    
\end{equation}
In fact, we will see that $t^*$ is the time at which the expected weight of the process (with densities of the reservoirs $p,q$) starting from the assumingly worst initial condition $\mathds{1}$ becomes ``close enough'' to the expected weight at equilibrium.
We prove that $t^*$ is a good estimate on the mixing time. \deleted{ and under the hypothesis of Theorem \ref{maintheorem}, $t^*$ is well approximated by $\dfrac{\log N}{2\pi^2}$.}
\begin{theorem}[Non asymptotic estimates]\label{maintheoremtwo}
    Under assumption \eqref{assumption}, there exists a universal constant $C$ such that for any $\epsilon \in ]0,1[$,  
    \begin{equation}\label{upperbound}
        t^* - C \left(1 + \log\left(\dfrac{1}{1 - \epsilon}\right)\right) \leq \tmix \leq t^* + C \left(1 + \log\left(\dfrac{1}{\epsilon}\right) + \log\left(\dfrac{1}{1-p}\right)\right).
    \end{equation}
\end{theorem}

% For numerous spin system, cutoff has been shown to occur at the time when the difference between the expected weight from the assumed worst initial configuration and the expected weight at equilibrium is smaller than the fluctuation of the weight at equilibrium: 
% \begin{equation*}
%     t^{**} = \minset{t \geq 0: \esperancewithstartingpoint{\mathds{1}}{S(\eta_t)} - \esperancewithstartingpoint{\pi}{S(\eta)} \leq \sqrt{\variancewithstartingpoint{\pi}{S(\eta)}}}.
% \end{equation*}
% Among many examples are the random walk on hypercube, see \cite{}, the Glauber dynamic of the Ising model, see \cite{Lubetzky2015, Lubetzky2016}. Our result is coherent with this, as we can prove that 
% \begin{equation*}
%     t^{**}-\dots \leq t^* \leq t^{**} + \dots 
% \end{equation*}

Estimating $t^*$ is a classical problem corresponding to the study of the discrete heat equation on the segment. It can be shown (see Appendix A in \cite{Goncalves2021} for example) that,

\begin{equation}\label{estimation_of_tmix} 
    t^* = \dfrac{1}{\pi^2} \log \left(\dfrac{N}{\sqrt{Np} \vee 1}\right) \pm \Ocal{1},
\end{equation}
 where the lower order term is bounded by some universal constant. This immediately implies the following corollary, of which Theorem \ref{maintheorem} is a direct consequence.
\begin{corollary}
    Under assumption \eqref{assumption}, the system exhibits cutoff at $t^*$ when 
    \begin{equation}\label{sufficient_condition_cutoff}
        \dfrac{1}{1-p} = N^{o(1)}.
    \end{equation}
    Moreover, if $p$ is bounded away from $1$, then the cutoff window is of order $\Ocal{1}$, and if $p$ is also bounded away from $0$, then $t^* = \dfrac{\log N}{2\pi^2} + \Ocal{1}$.  
\end{corollary}

We conjecture that cutoff should not be restrained by the condition \eqref{sufficient_condition_cutoff}.

\begin{conjecture}
    The system still exhibits cutoff at time $\dfrac{\log N}{2\pi^2}$ if $p = 1,\; q = 0$.
\end{conjecture}

\paragraph{Structure of the proofs. }
In Subsection $2.1$, we recall the negative dependence property of \orangeauthor{the} SSEP with reservoirs and the exponential bound introduced by Miller and Peres. In Subsection $2.2$, we \orangeauthor{collect} some elementary but useful estimates. In Subsection $2.3$, we introduce our method.  Section $3$ and Section $4$ are devoted to the proofs of the upper and lower bound in Theorem \ref{maintheoremtwo}, respectively. Finally, in Section $5$, we compute $t^*$ explicitly. 

\paragraph{Acknowledgment. }  The author warmly thanks Justin Salez for numerous fruitful discussions, thorough reading, and valuable comments on the draft. The author also thanks Hubert Lacoin for helpful discussions.

\section{Preliminary}
\subsection{Negative dependence property}
We recall the notion of negative dependence, which is essential throughout our proof.

\begin{definition}[Negative dependence]
    A random vector $Z = (Z_1, \dots, Z_n)$ taking values in $\{0,1\}^n$ is said to be negatively dependent (ND) if it satisfies
    \begin{equation}\label{negative_dependence}
        \forall A \subset [n], \;\; \esperance{\prod_{i \in A}Z_i} \leq \prod_{i \in A} \esperance{Z_i}.
    \end{equation}
\end{definition}

An important property of the SSEP with reservoirs is that it preserves the negative dependence property.
\begin{proposition}[The negative dependence property is preserved by SSEP with reservoirs, Lemma 12 in \cite{Salez2022}]\label{strongly_rayleigh_preserved_by_exclusion}
    Let the generator $\Lcal$ on $\Omega$ be defined by 
    \[\Lcal = \sum_{i,j = 1}^N a_{i, j} L_{i, j} + \sum_{i = 1}^N a^0_i L^0_i + \sum_{i = 1}^N a^1_i L^1_i,\]
    where 
    \begin{align*}
        L_{i,j}f(\eta ) &= f(\eta^{i \leftrightarrow j}) - f(\eta),\\
        L^0_i f(\eta) &= f(\eta^{i,0}) - f(\eta),\\ 
        L^1_i f(\eta) &= f(\eta^{i,1}) - f(\eta),
    \end{align*}
    where $a_{i, j}, a^0_i, a^1_i \in \Rbb_+, \, 1\leq i, j\leq N$. Then $\Lcal$ preserves the negative dependence property, \ie if $Z(0) \sim \mu$ is a negatively dependent vector for some measure $\mu$ on $\Omega$, then $Z(t) \sim \mu e^{\Lcal t}$ is also ND. 
\end{proposition}
In fact, the SSEP with reservoirs preserves a much richer property called the Strongly Rayleigh property, of which negative dependence is a consequence. We refer the readers to the beautiful paper $\cite{Borcea2009}$ for more details.

A particular case of inequality \eqref{negative_dependence} is when $|A| = 2$, which implies that the coordinates of an ND vector $Z$ are negatively correlated, and as a consequence, the weight of vector $Z$ is concentrated around its mean.
 
\begin{lemma}[Concentration of the weight of an ND vector]
    Let $Z = (Z_1, \dots, Z_n)$ be an ND vector. Let $S = \sum_{i = 1}^n Z_i$. Then 
    \begin{equation}\label{sum_of_strongly_Rayleigh_vector}
        \variance{S} \leq \esperance{S}.
    \end{equation}
\end{lemma}
% \begin{proof}
%     We simply write 
%     \begin{equation*}
%         \variance{S} = \sum_{i = 1}^n \sum_{j = 1}^n \covariance{Z_i}{Z_j}.
%     \end{equation*}
%     The conclusion follows from the following facts:
%     \begin{equation*}
%         \forall i \in [n], \; \covariance{Z_i}{Z_i} = \variance{Z_i} = \esperance{Z_i} - \esperance{Z_i}^2 \leq \esperance{Z_i}, 
%     \end{equation*}
%     and 
%     \begin{equation*}
%         \forall i \neq j,\; \covariance{Z_i}{Z_j} = \esperance{Z_iZ_j} - \esperance{Z_i}\esperance{Z_j} \leq 0,
%     \end{equation*}
%     by the negative dependence property.
% \end{proof}

\subsection{Some elementary estimates}
We first give a lemma about the symmetric simple random walk on the segment.
\begin{lemma}[Simple random walk on the segment]\label{one_walk_marginal}
    Let $(X(t))_{t\geq 0}$ be a continuous-time symmetric simple random walk on $\{0, \dots, N+1\}$, \orangeauthor{which jumps to the left (or to the right) at rate $N^2$}. For any $i \in \{0, \dots, N+1\}$, let $T_i$ be the first time that the walk reaches $i$, \ie $T_i = \minset{t \geq 0: X(t) = i}$. Then there exists a constant $c > 0$ (independent of $N$) such that for any $i \in \{0, \dots, N+1\}$, 
    \begin{align}
        \probawithstartingpoint{i}{T_0 \geq 2} &< e^{-c}. \label{time_to_reach_zero}
    \end{align}
\end{lemma}
\noindent This lemma can be proved by a classical hitting time estimate for non-negative supermartingale (see, \eg Proposition $2.1$ in \cite{Levin2017}, for the discrete version).\\

We write $\ber{p}$ for the Bernoulli distribution of parameter $p$. Now we recall a result about perturbation of product measures, first introduced by Miller and Peres in \cite{Miller2012} for the product of $\ber{1/2}$, extensively used by Lubetzky and Sly to prove cutoff for Ising model in a series of impressive papers \cite{Lubetzky2015}, \cite{Lubetzky2016}, \cite{Lubetzky2017}, and extended to the case of product of $\ber{p}$ by Salez in \cite{Salez2022}.
\begin{lemma}[Perturbation of the product lemma]\label{pertubation_product}
    Let $\Omega = \{0,1\}^n$. For each subset $S \subset [n]$, let $\phi_S$ be a distribution on $\{0,1\}^S$. Let $p \in (0,1)$, and let $\nu$ be the product measure $\ber{p}^{\otimes n}$ on $\Omega$. Let $\mu$ be the measure on $\Omega$ obtained by first sampling a subset $S \subset [n]$ via some measure $\mutilde$, and then, conditionally on $S$, generating independently the values on $S$ via $\phi_S$ and the values on $[n] \setminus S$ via $\ber{p}^{\otimes [n] \setminus S}$. Then 
    \begin{equation*}
        4 \dtv{\mu}{\nu}^2 \leq \norm{\dfrac{\mu}{\nu}-1}^2_{L^2(\nu)} \leq \esperance{\exponentialconstant^{\absolutevalue{S \cap S'}}} -1,
    \end{equation*}
    where $S, S'$ are \iid with law $\mutilde$, and $\exponentialconstant = \maxset{\dfrac{1}{p}, \dfrac{1}{1-p}}$.
\end{lemma}
We identify a subset $S \subset [n]$ with the vector $(\indicator{i \in S})_{1 \leq i \leq n}$. We remark here that the negative dependence property comes in very handy, as it allows us to bound the exponential moment in the above lemma by some quantity that depends only on the marginal of the random vector $S$, as stated in the following lemma.   
\begin{corollary}[Negative dependent perturbation of a product measure]\label{negative_dependent_pertubation}
    Under the above notations, if the random set $S$ is ND, then 
    \begin{equation*}
        4\dtv{\mu}{\nu}^2 \leq e^{(a-1)\sum_{i = 1}^n \proba{i \in S}^2} - 1. 
    \end{equation*}
\end{corollary}
For the proofs of Lemma \ref{pertubation_product} and Corollary \ref{negative_dependent_pertubation}, see Lemma $9$ in \cite{Salez2022}.
\subsection{Framework and some definitions}\label{framework}    
The SSEP with reservoirs $(\eta_t)_{t\geq 0}$ evolves according to the following transitions: 
\begin{enumerate}
    \item $\eta \mapsto \eta^{i \leftrightarrow i + 1}$ (exchange between site $i$ and site $i+1$), which occurs at rate $N^2$,
    \item Resampling the value at site $1$ by an independent Bernoulli $\ber{p}$, which occurs at rate $N^2$, 
    \item Resampling the value at site $N$ by an independent Bernoulli $\ber{q}$, which occurs at rate $N^2$. 
\end{enumerate}
We introduce another Markov process closely related to the SSEP with reservoirs.
\paragraph{The colored interchange process.}

Let $\setcolorinterchange := \symmetricgroup \times \{R,B,G\}^N$, where $\symmetricgroup$ is the symmetric group on $[N]$, and $R,B,G$ stand for red, blue, and green. Each element $(\sigma, b) \in \setcolorinterchange$ describes a way to put $N$ colored labelled individuals on the segment as follows. The individuals are labelled $1,2, \dots, N$. For any $i \in [N]$, the individual labelled $i$ is located at site $\sigma(i)$ and is colored $b(i)$. \textit{The colored interchange process} X is defined as the Markov process taking values in the state space $\setcolorinterchange$ which evolves according to the following transition: 
\begin{enumerate}
    \item $(\sigma, b) \mapsto \left((i, i+1) \circ \sigma,b\right)$ (the individuals at sites $i$ and $i+1$ exchange their positions), which occurs at rate $N^2$, for $1\leq i \leq N-1$,
    \item $(\sigma, b) \mapsto (\sigma, b^{\sigma^{-1}(1),B})$ (recoloring the individual at site $1$ blue), which occurs at rate $N^2$, 
    \item $(\sigma, b) \mapsto (\sigma, b^{\sigma^{-1}(N),G})$ (recoloring the individual at site $N$ green), which occurs at rate $N^2$.
\end{enumerate}
There is a natural coupling between $(\eta_t)_{t\geq 0}$ and $X$ as follows.
\paragraph{Natural coupling between $(\eta_t)_{t\geq 0}$ and $X$.} A coupling of the two processes is given by making the transitions $1, 2, 3$ listed above of the two processes $(\eta_t)_{t\geq 0}$ and $X$ occur at the same time.\\
\hfill\\
Roughly speaking, the labels and colors are added to keep better track of the exchange of information inside the bulk and to memorize which reservoirs the resamplings come from. For an introduction to the interchange process and its relation with the exclusion process, see \cite{Levin2017}, chapter $23$.

% \paragraph{The "pushforward" measure.} For $\zeta \in \Omega$ and $x = (\sigma, b) \in \setcolorinterchange$, let $\mu(\zeta, x)$ be the product law on $\Omega$ given by:
% \begin{equation}\label{pushforward_measure}
%     \esperancewithstartingpoint{\mu(\zeta, x)}{\eta(i)} = \begin{cases}
%     \zeta(\sigma(i)) &\text{if } b_i = R,\\
%     p &\text{if } b_i = B,\\
%     q &\text{if } b_i = G.
%     \end{cases}
% \end{equation}

\paragraph{The ``pushforward" function} Let the function $\pushforwardfunction: \Omega \times \setcolorinterchange \times \Omega \times \Omega \to \Omega$ be defined by, for any $x = (\sigma, b) \in \setcolorinterchange$ and $\eta, v^B, v^G \in \Omega$, 
\begin{equation}
    \pushforwardfunction(\eta,x,v^B,v^G)(i) = \begin{cases}
    \eta(\sigma^{-1}(i)) &\text{if } b\left(\sigma^{-1}(i)\right) = R,\\
    v^B(i) &\text{if } b\left(\sigma^{-1}(i)\right) = B,\\
    v^G(i) &\text{if } b\left(\sigma^{-1}(i)\right) = G.
    \end{cases}
\end{equation}

The interest of introducing the process $X$ is the following lemma, whose proof is straightforward from the definition of the natural coupling.
\begin{lemma}\label{natural_coupling}
    Let $x_0 \in \setcolorinterchange$ be the configuration where for any $i\in [N]$, the individual labelled $i$ is located at site $i$ and is colored red:
    \[x_0 = \left(\identity, \redvector \right),\]
    Let $X, \xi^B, \xi^G$ be independent and as follows.
    \begin{itemize}
        \item X is a colored interchange process starting at $x_0$.
        \item $\xi^B \sim \ber{p}^{\otimes N}$.
        \item $\xi^G \sim \ber{q}^{\otimes N}$.
    \end{itemize} 
    Let $(\eta_t)_{t\geq 0}$ be the SSEP with reservoirs started from some configuration $\eta \in \Omega$. Then for any $t \geq 0$, 
    \begin{equation}
        \eta_t \overset{d}{=} \pushforwardfunction(\eta, X_t, \xi^B, \xi^G), 
    \end{equation}
    where $\overset{d}{=}$ means equal in distribution.
\end{lemma}

\paragraph{Red, blue, and green regions.} For any $x \in \setcolorinterchange$, we denote by $R(x)$ \textbf{the red region}, \ie the set of the sites containing the red individuals:
\begin{equation*}
    R(x):= \{i \in [N] \big| b(\sigma^{-1}(i)) = R\}.
\end{equation*}
The blue region $B(x)$ and the green region $G(x)$ are defined similarly. As the two reservoirs recolor the particles blue or green, the red region evolves exactly as a SSEP with two reservoirs of density $0$. Then $t^*$ is the time at which the red region becomes small enough:
\begin{equation}\label{t*andredregion}
    t^* = \minset{t \geq 0: \esperancewithstartingpoint{x_0}{\absolutevalue{R(X_t)}} \leq \sqrt{Np} \vee 1},   
\end{equation}
with $x_0$ as in Lemma \ref{natural_coupling}. We present a graphical construction of $X$ that allows us to reveal the green region before the red and blue regions.
\paragraph{Graphical construction of the colored interchange process.} We can construct $X$ in the following way.
\begin{equation}
    X = \Psi (\sigma, b, \xisamplingone, \xisamplingn, , (\Xi^G_{i})_{1 \leq i \leq N-1}, (\Xi^{BR}_{i})_{1 \leq i \leq N- 1}),
\end{equation}
where $(\sigma, b) \in \setcolorinterchange$ and $\xisamplingone, \xisamplingn, (\Xi^G_{i})_{1 \leq i \leq N-1}, (\Xi^{BR}_{i})_{1 \leq i \leq N- 1}$ are independent and as follows.
\begin{itemize}
    \item $\xisamplingone$ and $\xisamplingn$ are homogeneous Poisson processes of intensity $N^2 \drm t$ which indicate the times at which we recolor the individuals at site $1$ and site $N$, respectively. 
    \item $\Xi^G_{i}$ and $\Xi^{BR}_{i}$ are homogeneous Poisson processes of intensity $N^2 \drm t $, $1 \leq i \leq N - 1$. Each time $\Xi^G_{i}$ jumps, the two individuals at sites $i$ and $i+1$ exchange their positions if at least one of them is green, and each time $\Xi^{BR}_{i}$ jumps, the two individuals at sites $i$ and $i+1$ exchange their positions if none of them is green.
\end{itemize}
This construction gives us the colored interchange process $X$ with initial condition $(\sigma, b)$.

\paragraph{Trajectories of a single individual. }It is well known, see \eg Chapter 23 of \cite{Levin2017}, that if we observe the trajectory of a single labelled individual, we see a continuous-time simple random walk on the segment where the conductance of any edge is $N^2$. As we add the colors here, we also see that the labelled individual is always recolored at rate $N^2$ when it is at site $1$ or $N$. 
\paragraph{Trajectories of the green regions. }
The above construction allows us to reveal $(G(X_t))_{t\geq 0}$ before $(R(X_t))_{t\geq 0}$ and $(B(X_t))_{t\geq 0}$. In fact, $(G(X_t))_{t\geq 0}$ is measurable with respect to the $\sigma$-algebra $\greentrajectory$ generated by $\xisamplingone, \xisamplingn$, and  $\Xi^G$.
 
 \paragraph{Number of crossings.} Let $(L_t)_{t\geq 0}$ be the process that counts the number of times that a blue or red individual is recolored green. Note that if an individual is recolored blue at site $1$, then it needs to cross the bulk to be recolored green. Accordingly, we call $(L_t)_{t\geq 0}$ the number of crossings. A simple but important observation is that $(L_t)_{t\geq 0}$ is also $\greentrajectory$-measurable.
 
 \section{The upper bound}
Since we do not have an explicit formula for the invariant measure, we will compare two processes from two arbitrary configurations $\eta$ and $\tilde{\eta}$ and use the fact that 
\begin{equation}
    d(t) \leq \max_{\eta, \tilde{\eta}\in \Omega} \dtv{\probawithstartingpoint{\eta}{\eta_t \in \cdot}}{\probawithstartingpoint{\tilde{\eta}}{\eta_t \in \cdot}},
\end{equation}
which is due to the convexity of the total variation distance. Our goal now is to compare the distributions of those two processes at our predicted mixing time $t^*$. Our strategy is to match perfectly the green region of the two processes by the graphical construction above and to view the distributions on the remaining sites as a product measure perturbed by the red region to compare them using Lemma \ref{pertubation_product}. 
% \paragraph{Number of label of a color. }We denote by $W^i(x)$ the number of labels of color $i$ in the configuration $x \in \setcolorinterchange$, where $i\in \{R, B, G\}$. 

% \hfill\\
We will need the following lemmas and propositions.
\begin{lemma}[Exponential decay of the red region]\label{exponential_decay_w0} 
    For any $x \in \setcolorinterchange$, for any $t \geq 0$, 
    \begin{equation*}
        \esperancewithstartingpoint{x}{\absolutevalue{R(X_{2t})}} \leq e^{-c\lfloor t \rfloor} \absolutevalue{R(x)},
    \end{equation*}
    where $c$ is the constant in Lemma \ref{one_walk_marginal}.
\end{lemma}
\begin{lemma}[Fast increase of the number of crossings]\label{lemma_number_of_crossing}
    There exists a constant $C$ such that for $t_2 = C(1 + \log(1/\epsilon))$, for any initial configuration $x \in \setcolorinterchange$, 
    \begin{equation*}
        \probawithstartingpoint{x}{ L_{t_2} < 2N} \leq \epsilon/4.
    \end{equation*}
\end{lemma}

\begin{proposition}[Negative dependence property of conditional law]\label{conditional_negative_dependence}
    For any initial configuration $x \in \setcolorinterchange$, almost surely, conditionally on $\greentrajectory$, $R(X_t)$ is negatively dependent at any time $t \geq 0$.
\end{proposition}

\begin{lemma}[Conditional anticoncentration inequality]\label{lemma_anticoncentration_inequality}
    Let $t$ be a positive number, and let $x = (\sigma, b) \in \setcolorinterchange$ be an initial configuration. For any individual $i$ that is colored blue or red in $x$, for any site $j$, on the event $\{L_t \geq 2N\}$, 
    \begin{equation*}
        \probawithstartingpoint{x}{\sigma_t(i) = j, \sigma_s(i) \notin G(X_s), \forall0 \leq s \leq t\big| \greentrajectory} \leq \dfrac{1}{N}.
    \end{equation*}
    % where the meaning of the event is that the individual $i$ ends up at site $j$ at time $t$ and was never recolored green before.
\end{lemma}
The following is a direct consequence, obtained by summing the inequality in Lemma \ref{lemma_anticoncentration_inequality} over all red individuals in configuration $x$. 
\begin{corollary}[Conditional marginal of $R(X_t)$]\label{marginal_red_conditional}
    Let $t$ be a positive number, and let $x \in \setcolorinterchange$ be an initial configuration. For any site $j$, on the event $\{L_t \geq 2N\}$, 
    \begin{equation*}
        \probawithstartingpoint{x}{j \in R(X_{t})\big| \greentrajectory} \leq \dfrac{\absolutevalue{R(x)}}{N}.
    \end{equation*}
\end{corollary}

Now we are ready to prove the upper bound in Theorem \ref{maintheoremtwo}.
\begin{proof}[Proof of the upper bound]
    Let $\eta, \tilde{\eta} \in \Omega$ arbitrary, and let $c$ be the constant in Lemma \ref{one_walk_marginal}.
    % We construct two processes $(\eta_t)_{t \geq 0}$ starting from $\eta_0 = \eta$ and $(\tilde{\eta}_t)_{t \geq 0}$ starting from $\tilde{\eta}_0 = \tilde{\eta} $ via the colored interchange process $X$ starting at $x_0 = (\textrm{id}, (R, \dots, R))$ (initially every site $i$ is given the individual $i$ and every individual is colored red).
    Let $x_0, X, \eta^B, \eta^G$ be defined as in Lemma \ref{natural_coupling}. Let $(\zeta_t)_{t\geq 0}$ and $(\zetatilde_t)_{t\geq 0}$ be defined by 
    \begin{align*}
        &\zeta_t = \pushforwardfunction(\eta, X_t, \xi^B, \xi^G), \\
        &\zetatilde_t = \pushforwardfunction(\tilde{\eta}, X_t, \xi^B, \xi^G).
    \end{align*}
    For $x \in \Omega$, we denote by $\probawithstartingpoint{x}{\cdot}$ the law of the process $X$ starting from $x$. By Lemma \ref{natural_coupling}, 
    \begin{align*}
        &\probawithstartingpoint{\eta}{\eta_t \in \cdot} = \probawithstartingpoint{x_0}{\zeta_t \in \cdot},\\
        &\probawithstartingpoint{\tilde{\eta}}{\eta_t \in \cdot} = \probawithstartingpoint{x_0}{\zetatilde_t \in \cdot}.
    \end{align*}
    So now we can compare the distributions of $\zeta_t$ and $\zetatilde_t$ instead of those of $\eta_t$ and $\tilde{\eta}_t$.
    We divide the proof into two cases: $Np \leq 1$ and $Np > 1$.
    \paragraph{Case 1: $Np \leq 1$.} We see that if $\absolutevalue{R(X_t)} = 0$, then $\probawithstartingpoint{x_0}{\zeta_t \in \cdot \big| X_t} = \probawithstartingpoint{x_0}{\tilde{\zeta_t} \in \cdot\big| X_t}$, by definition of $\pushforwardfunction$. Hence for $t = t^* + 2m$, for some $m \in \Zbb_+$,
    \begin{align*}
        &\dtv{\probawithstartingpoint{x_0}{\zeta_t \in \cdot}}{\probawithstartingpoint{x_0}{\zetatilde_t \in \cdot}}\\
        &\leq \esperance{\norm{\probawithstartingpoint{x_0}{\zeta_t \in \cdot \big| X_t} - \probawithstartingpoint{x_0}{\tilde{\zeta}_t \in \cdot \big|X_t}}_{TV}}\\
        &\leq \probawithstartingpoint{x_0}{\absolutevalue{R(X_t)} > 0} \\
        &\leq \esperancewithstartingpoint{x_0}{\absolutevalue{R(X_t)}} \leq e^{-cm} \esperancewithstartingpoint{x_0}{\absolutevalue{R(X_{t^*})}} = e^{-c m},    
    \end{align*}
    where the first inequality is due to Jensen's inequality, the second inequality is by upper bounding $\norm{\cdot}_{TV}$ by $1$ on the event $\{R(X_t) > 0\}$, the third inequality is because $\absolutevalue{R(X_t)} \in \Zbb_+$, the last inequality is by Lemma \ref{exponential_decay_w0}, and the equality is by \eqref{t*andredregion}. We can take $m = \left\lceil\dfrac{\log \dfrac{1}{\epsilon}}{c}\right\rceil$ to make $e^{-cm}$ smaller than $\epsilon$, which finishes the proof.
    \paragraph{Case 2: $Np > 1$.}
    Let $t_1, t_2, \alpha$ be some positive numbers that we will choose later. We see that,
    \begin{align}\label{markov_at_time_t1}
        &\dtv{\probawithstartingpoint{x_0}{\zeta_t \in \cdot }}{\probawithstartingpoint{x_0}{\zetatilde_t \in \cdot }} \nonumber\\
        &\leq \esperancewithstartingpoint{x_0}{\norm{\probawithstartingpoint{x_0}{\zeta_t \in \cdot |X_{t_1}} - \probawithstartingpoint{x_0}{\zetatilde_t \in \cdot |X_{t_1}}}_{TV}} \nonumber\\
        &= \esperancewithstartingpoint{x_0}{\norm{\probawithstartingpoint{X_{t_1}}{\zeta_{t_2} \in \cdot} -\probawithstartingpoint{X_{t_1}}{\zetatilde_{t_2} \in \cdot}  }_{TV}} \nonumber\\
        &\leq \probawithstartingpoint{x_0}{\absolutevalue{R(X_{t_1})} > \alpha} + \max_{x: \absolutevalue{R(x)}\leq \alpha} \norm{\probawithstartingpoint{x}{\zeta_{t_2} \in \cdot} - \probawithstartingpoint{x}{\zetatilde_{t_2} \in \cdot}}_{TV}, 
    \end{align}
    where the first inequality is by Jensen's inequality, the equality is due to the Markov property of the process $X$ at time $t_1$, and the second inequality is by upper bounding the total variation distance by $1$ on the event $\{\absolutevalue{R(X_{t_1})} > \alpha\}$. For any $x \in \setcolorinterchange$, we use the graphical construction in Section $2$ to construct the process $X$ starting from $x$. Then
    \begin{align}\label{conditionally_on_greentrajectory}
        &\norm{\probawithstartingpoint{x}{\zeta_{t_2} \in \cdot} - \probawithstartingpoint{x}{\zetatilde_{t_2} \in \cdot}}_{TV} \nonumber\\
        &\leq \esperance{\norm{\probawithstartingpoint{x}{\zeta_{t_2} \in \cdot|\greentrajectory} - \probawithstartingpoint{x}{\zetatilde_{t_2} \in \cdot|\greentrajectory}}_{TV}} \nonumber\\
        &\leq \probawithstartingpoint{x}{L_{t_2} < 2N} +  \esperance{\norm{\probawithstartingpoint{x}{\zeta_{t_2} \in \cdot|\greentrajectory} - \probawithstartingpoint{x}{\zetatilde_{t_2} \in \cdot|\greentrajectory}}_{TV} \indicator{ L_{t_2} \geq 2N}},
    \end{align}
    where the first inequality is by Jensen's inequality, and the second inequality is simply by upper bounding the total variation distance by $1$ on the event $\{L_{t_2} < 2N\}$. Combining \eqref{markov_at_time_t1} and \eqref{conditionally_on_greentrajectory}, we deduce that $\dtv{\probawithstartingpoint{x_0}{\zeta_t \in \cdot }}{\probawithstartingpoint{x_0}{\zetatilde_t \in \cdot }}$ does not exceed 
    \begin{align}\label{bounded_by_the_sum}
        &\probawithstartingpoint{x_0}{\absolutevalue{R(X_{t_1})} > \alpha} + \max_{x \in \setcolorinterchange} \probawithstartingpoint{x}{L_{t_2} < 2N} \nonumber\\
        &\;\; +\max_{x: \absolutevalue{R(x)} \leq \alpha} \esperance{\norm{\probawithstartingpoint{x}{\zeta_{t_2} \in \cdot|\greentrajectory} - \probawithstartingpoint{x}{\zetatilde_{t_2} \in \cdot|\greentrajectory}}_{TV} \indicator{ L_{t_2} \geq 2N}}.
    \end{align}
    We separately estimate the three terms in the sum \eqref{bounded_by_the_sum}.\\
    \hfill\\
    \textbf{The first term.} We choose $t_1 = t^* + 2m$ for some $m \in \Zbb_+$ that we will choose later. Then by Lemma \ref{exponential_decay_w0} and equality \eqref{t*andredregion}, $\esperancewithstartingpoint{x_0}{\absolutevalue{R(X_{t_1})}} \leq \esperancewithstartingpoint{x_0}{\absolutevalue{R(X_{t^*})}} e^{-cm} = \sqrt{Np} e^{-cm}$. Hence by Markov's inequality, 
    \[\probawithstartingpoint{x_0}{\absolutevalue{R(X_{t_1})} \geq \alpha } \leq \dfrac{\sqrt{Np}e^{-cm}}{\alpha}.\]
     \textbf{The second term.} This term is smaller than $\epsilon/4$ for $t_2$ as in Lemma \ref{lemma_number_of_crossing}.\\
    \hfill\\
    \textbf{The third term.} Let $x\in \setcolorinterchange$. Observe that, under $\Pbb_x$, $R(X_{t_2}) \cup B(X_{t_2}) = [N] \setminus G(X_{t_2})$ is $\greentrajectory$-measurable. We write $\zetaredblue$ for $\zeta_{t_2}(\redbluettwo)$ and $\zetaredbluetilde$ for $\zetatilde_{t_2}(\redbluettwo)$. Conditionally on $\greentrajectory$, by construction, the distribution of $\zeta_{t_2}$ and $\zetatilde_{t_2}$ on $G(X_{t_2})$ is a product of $\ber{q}$, independent of the restriction of $\zeta_{t_2}$ and $\zetatilde_{t_2}$ on $\redbluettwo$, so we can safely project onto $\redbluettwo$ to obtain 
    \begin{align}\label{conditional_projection}
        &\esperance{\norm{\probawithstartingpoint{x}{\zeta_{t_2} \in \cdot|\greentrajectory} - \probawithstartingpoint{x}{\zetatilde_{t_2} \in \cdot|\greentrajectory}}_{TV} \indicator{ L_{t_2} \geq 2N}} \nonumber\\
        &= \esperance{\norm{\probawithstartingpoint{x}{\zetaredblue \in \cdot|\greentrajectory} - \probawithstartingpoint{x}{\zetaredbluetilde \in \cdot|\greentrajectory}}_{TV}\indicator{ L_{t_2} \geq 2N}} \\
        &\leq 2 \sup_{\eta} \esperance{\norm{\probawithstartingpoint{x}{\zetaredblue \in \cdot|\greentrajectory} - \nu_{\redbluettwo}}_{TV}\indicator{ L_{t_2} \geq 2N}}, \nonumber
    \end{align}
    where $\nu_{\redbluettwo}$ is the product measure of $\ber{p}$ on $\redbluettwo$, by the triangle inequality. By construction of $\zeta$, conditionally on $\greentrajectory$, the variable $\zetaredblue$ can be constructed by first sampling the set $R(X_{t_2})$, then sampling the values of $\zeta_{t_2}$ on $R(X_{t_2})$ conditionally on $R(X_{t_2})$, and then sampling independently the values of $\zeta_{t_2}$ on $B(X_{t_2})$ by a product of $\ber{p}$. So thanks to Corollary \ref{negative_dependent_pertubation}, for $\exponentialconstant = \dfrac{1}{\min\{p, 1-p\}}$,
    \begin{equation}\label{application_of_negative_pertubation}
    \begin{split}
        & 2\norm{\probawithstartingpoint{x}{\zetaredblue \in \cdot|\greentrajectory} - \nu_{\redbluettwo}}_{TV}\\
        & \leq \sqrt{  \exp \left(\sum_{i \in \redbluettwo}(\exponentialconstant - 1) \probawithstartingpoint{x}{i \in R(X_{t_2})\big| \greentrajectory}^2\right) - 1},
    \end{split}
    \end{equation}
    By Corollary \ref{marginal_red_conditional}, on the event $\{L_{t_2} \geq 2N\}$,
    \begin{equation}\label{sum_of_proba_to_be_red}
        \sum_{i \in \redbluettwo} \probawithstartingpoint{x}{i \in R(X_{t_2})\big| \greentrajectory}^2 \leq \absolutevalue{\redbluettwo} \dfrac{\absolutevalue{R(x)}^2}{N^2} \leq N \dfrac{\absolutevalue{R(x)}^2}{N^2} = \dfrac{\absolutevalue{R(x)}^2}{N}.
    \end{equation}
    \eqref{conditional_projection}, \eqref{application_of_negative_pertubation}, and \eqref{sum_of_proba_to_be_red} together imply that the third term is upper bounded by $\sqrt{e^{\dfrac{(\exponentialconstant - 1)\alpha^2}{N}} - 1}$.\\
    So altogether, with $t_1 = t^* + 2m$ and $t_2$ as in Lemma \ref{lemma_number_of_crossing}, 
    \begin{equation}\label{bounded_by_the_sum_two}
        \dtv{\probawithstartingpoint{x_0}{\zeta_t \in \cdot }}{\probawithstartingpoint{x_0}{\zetatilde_t \in \cdot }}  \leq \dfrac{\sqrt{Np}e^{-cm}}{\alpha} + \dfrac{\epsilon}{4} + \sqrt{e^{\dfrac{(\exponentialconstant - 1)\alpha^2}{N}} - 1}.
    \end{equation}
    We take $\alpha = \dfrac{\epsilon\sqrt{N}}{\sqrt{8\exponentialconstant}}$. Then
    \[\dfrac{\exponentialconstant - 1}{N} \alpha^2 \leq \dfrac{\epsilon^2}{8} \leq \log \left(1 + \dfrac{\epsilon^2}{4}\right),\] where we have used the inequality $\log(1+x) \geq x/2$ if $0 < x < 1$. This implies that the last term on the right-hand side of \eqref{bounded_by_the_sum_two} is smaller than $\epsilon/2.$ The first term on the right-hand side of \eqref{bounded_by_the_sum_two} now becomes $\dfrac{\sqrt{8 \exponentialconstant p}}{\epsilon} e^{-cm}$. Observe that
    \begin{equation*}
        ap = p \maxset{\dfrac{1}{p}, \dfrac{1}{1-p}} = \maxset{1, \dfrac{p}{1-p}} \leq \dfrac{1}{1-p},
    \end{equation*}
    Hence
    \begin{equation*}
        \dfrac{\sqrt{8 \exponentialconstant p}}{\epsilon} e^{-cm} \leq\dfrac{1}{\epsilon} \sqrt{\dfrac{8}{1-p}} e^{-cm}.    
    \end{equation*}
    We can take 
    \[m = \left\lceil\dfrac{1}{c}\left(\dfrac{1}{2} \log (128) + \dfrac{1}{2} \log \left(\dfrac{1}{1-p}\right) + 2 \log(1/\epsilon)\right)\right\rceil\]
    to make that term smaller than $\epsilon/4$. Hence for $t = t^* + 2m + t_2 $ with $m$ and $t_2$ defined as above, $\dtv{\probawithstartingpoint{\eta}{\eta(t) \in \cdot}}{\probawithstartingpoint{\tilde{\eta}}{\tilde{\eta}(t) \in \cdot}} < \epsilon$, for any $\eta, \tilde{\eta} \in \Omega$, which finishes our proof.
\end{proof}

The rest of this section is dedicated to proving Lemma \ref{exponential_decay_w0}, Lemma \ref{lemma_number_of_crossing}, Lemma \ref{lemma_anticoncentration_inequality}, and Proposition \ref{conditional_negative_dependence}. First, we prove Lemma \ref{exponential_decay_w0}.
\begin{proof}[Proof of Lemma \ref{exponential_decay_w0}]
    As we mentioned near the end of Subsection \ref{framework}, each labelled individual moves as a continuous-time simple random walk on the segment where every edge has conductance $N^2$. Each individual is also recolored at site $1$ and $N$ at rate $N^2$. Hence the time at which a red individual is recolored is also the time that a continuous-time random walk on $\Zbb$, whose edges are given conductance $N^2$, starting from the same site reaches $0$ or $N+1$ (we imagine that the red individual jumps to site $0$ or site $N+1$ when it is recolored). Hence by Lemma \ref{one_walk_marginal}, the probability that an individual remains red up to time $2$ is smaller than $e^{-c}$. Summing over all red individuals, we get 
    \[\esperancewithstartingpoint{x}{\absolutevalue{R(X_2)}} \leq e^{-c} \absolutevalue{R(x)}.\]
    Moreover, since $\absolutevalue{R(X_t)}$ stochastically decreases with respect to $t$ as the individuals are only recolored blue or green,
    the conclusion is obtained simply by iterating the above inequality via the Markov property.
\end{proof}
To prove Lemma \ref{lemma_number_of_crossing}, Proposition \ref{conditional_negative_dependence}, and Lemma \ref{lemma_anticoncentration_inequality}, we need to understand the evolution of the system conditionally on $\greentrajectory$. We start by describing how the system evolves once we fix a realization of $\xisamplingone, \xisamplingn$, and $\Xi^G$.
\begin{observation}[Evolution of the red and blue individuals conditionally on $\greentrajectory$]\label{observation_conditional_evolution}
Let $T_0 = 0$, and let $(T_i)_{i \geq 1}$ be the times at which a point of $\xisamplingone, \xisamplingn$, or $\Xi^G$ appears, which are $\greentrajectory$-measurable and strictly increase to infinity almost surely. From the graphical construction of $X$, we see that conditionally on $\greentrajectory$, when $\Xi^{BR}$ is revealed, the red and blue individuals evolve as follows. On any time interval $(T_i, T_{i+1}), \, i \geq 0,$ two neighbor (red or blue) individuals exchange their positions at rate $N^2$, and at time $T_i, \, i\geq 1$, the system is forced to take some transitions by the environment, as follows. 
\begin{itemize}
    \item Two neighbor green individuals exchange their positions. Then nothing happens to the red and blue individuals,
    \item A red or blue individual, say at site $i$, is forced to exchange positions with a green individual, say at site $i+1$, due to a point of $\Xi^G$. 
    \item The reservoirs recolor a green individual blue at site $1$ or recolor a blue (or red) individual green at site $N$, due to a point of $\xisamplingone$ or $\xisamplingn$.
\end{itemize}
In short, the blue and red individuals evolve as a simple exclusion process conditionally on the environment created by the green region.
\end{observation}
Now we prove Proposition \ref{conditional_negative_dependence}. 
\begin{proof}[Proof of Proposition \ref{conditional_negative_dependence}]
    The proof is inspired by that of Proposition \ref{strongly_rayleigh_preserved_by_exclusion}. We fix a realization of $\xisamplingone, \xisamplingn, \Xi^G$ and let $(T_i)_{i\geq 0}$ be defined as in Observation \ref{observation_conditional_evolution}. We will prove the following two statements, which are conditional on $\greentrajectory$:
    \begin{enumerate}
        \item For any $i \geq 0$, if $R(X_{T_i})$ is ND, then $R(X_t)$ is also ND, for any $t \in [T_i, T_{i+1}[$.
        \item For any $i \geq 0$, if $R(X_{T_{i}-})$ is ND, then $R(X_{T_i})$ is also ND.
    \end{enumerate}
    These two statements clearly imply that if $R(X_0)$ is ND (which is the case when $X_0$ is deterministic), then at any time $t \geq 0$, $R(X_t)$ is still ND. Now we prove the first statement. 
    % On $[T_i, T_{i+1}[$, when we realize $\Xi^{BR}$, the red labels evolve as an exclusion process conditionally on avoiding $G(X_{T_i})$.
    %Here we want to use another graphical construction. Instead of putting the Poisson clocks on the particles, we will put Poisson clocks of intensity $N^2$ on the edges of the graph that do not touch the particles of type $2$. Each time a clock on an edge ring, we will exchange the particles at the two endpoints of that edge. This graphical construction gives the same evolution (in distribution) of the particles of type $0$ as our process at the beginning. 
    From Observation \ref{observation_conditional_evolution}, we see that on the time interval $[T_i, T_{i+1}[$, the red region evolves according to the generator $\tilde{L}$ given by
    \[\tilde{L} = N^2\sum_{i} L_{i, i+1},\]
    with $L_{i, i+1}$ as in Proposition \ref{strongly_rayleigh_preserved_by_exclusion}, where the sum is taken over all the site $i$ such that $\{i, i+1\} \cap G(X_{T_i}) = \emptyset$. We can now apply Proposition \ref{strongly_rayleigh_preserved_by_exclusion} to conclude that the ND property of the red region is preserved from time $T_i$ to $T_{i+1}-$. Now we prove the second statement. According to Observation \ref{observation_conditional_evolution}, at $T_i$, only a few following things can happen.
    \begin{itemize}
        \item Two green individuals exchange their positions. Then nothing happens to the red region, \ie $R(X_{T_i}) = R(X_{T_i-})$. It is straightforward that ND property is preserved.
        \item A red or blue individual is forced to exchange positions with a green individual, say the exchange happens between two sites $j$ and $j+1$. Then $R(X_{T_i}) = R(X_{T_i-})^{j \leftrightarrow j+1}$. This does not affect the ND property as the inequality \eqref{negative_dependence} for $R(X_{T_i})$ is exactly that for $R(X_{T_i-})$ when we replace $A$ by $A^{j \leftrightarrow j+1}$.
        \item The reservoirs recolor an individual, for example, at site $1$. This corresponds to setting the first coordinate to $0$. If $1 \notin A$, this operation does not affect the inequality \eqref{negative_dependence}. If $1 \in A$, the two sides of inequality \eqref{negative_dependence} become $0$. In both cases, this operation preserves the inequality \eqref{negative_dependence} and hence the ND property.
    \end{itemize}
    This finishes our proof.
\end{proof}

\begin{proof}[Proof of Lemma \ref{lemma_number_of_crossing}]
    Let $(T_i)_{i\geq 0}$ be as in Observation \ref{observation_conditional_evolution}. Let $m \in \Nbb$ be a number that we will choose later. We call \textit{the modified dynamics} the evolution of the system where we close the reservoir at site $N$ during the time interval $[0,2m]$ and close the reservoir at site $1$ during time $]2m, 4m]$, \ie we ignore the points of $\xisamplingn$ on the time interval $[0, 2m]$ and the points of $\xisamplingone$ on $]2m, 4m]$. We claim that at time $4m$, we have fewer blue and red individuals recolored green in the modified dynamics than in the original dynamics (actually, closing some reservoirs during any time only decreases $(L_t)_{t\geq 0}$). More precisely, suppose that in the modified dynamics, some individual labelled $i$ is recolored blue while being green at some time $T_k$ (or simply a blue or red individual if $k = 0$) and then keeps its color until being recolored green at some time $T_l$ for some $l > k$. Then in the original dynamics, the individual labelled $i$ is still green at time $T_k-$ and is recolored blue at time $T_k$ (or simply blue or red if $k = 0$), and after that either it is recolored green at some time in the interval $[0,2m]$, or it stays blue (or red) until time $2m$ and then being recolored green at time $T_l$, which proves the observation. Now we prove that in the modified dynamics, there are many recolorings occurring at site $N$. The point here is that when we close the reservoir at site $N$, the other reservoir at site $1$ quickly paints almost the whole bulk blue, and vice versa, when we close the reservoir at site $1$, the reservoir at site $N$ quickly paints almost the whole bulk green, and hence there are many blue individuals recolored green. More precisely, we write $\esperancewithstartingpointmodified{}{\cdot}, \probawithstartingpointmodified{}{\cdot}, \variancewithstartingpointmodified{}{\cdot}$ for the expectation, the probability, and the variance taken with respect to the modified dynamics. In the modified dynamics, for any initial configuration $x$, the set $G(X_{2m})$ is ND, by the same argument as in Proposition \ref{conditional_negative_dependence}. Hence by Lemma \ref{strongly_rayleigh_preserved_by_exclusion},
    \[\variancewithstartingpointmodified{x}{\absolutevalue{G(X_{2m})}} \leq \esperancewithstartingpointmodified{x}{\absolutevalue{G(X_{2m})}}.\]
    By Lemma \ref{one_walk_marginal} and a proof similar to that of Lemma \ref{exponential_decay_w0}, for any $m \in \Zbb_+$,
    \[\esperancewithstartingpointmodified{x}{\absolutevalue{G(X_{2m})}}\leq \absolutevalue{G(x)} e^{-cm} \leq Ne^{-cm}.\]
    Then for $m =\left\lceil \dfrac{1}{c} \log \dfrac{1000}{\epsilon}\right\rceil$, 
    \begin{align*}
        \probawithstartingpointmodified{x}{\absolutevalue{G(X_{2m})} \geq N/4} &= \probawithstartingpointmodified{x}{\absolutevalue{G(X_{2m})} - \esperancewithstartingpointmodified{x}{\absolutevalue{G(X_{2m})}} \geq N/4 - \esperancewithstartingpointmodified{x}{\absolutevalue{G(X_{2m})}}} \\
        &\leq \probawithstartingpointmodified{x}{\absolutevalue{G(X_{2m})} - \esperancewithstartingpointmodified{x}{\absolutevalue{G(X_{2m})}} \geq N/4 - Ne^{-cm}}\\
        &\leq \dfrac{\variancewithstartingpointmodified{x}{\absolutevalue{G(X_{2m})}}}{(1/4 - e^{-cm})^2N^2} \leq \dfrac{Ne^{-cm}}{(1/4 - e^{-cm})^2N^2} \leq \dfrac{\epsilon}{32}.             
    \end{align*}
    By the same argument, 
    \[\probawithstartingpointmodified{x}{\absolutevalue{G(X_{4m})} \leq 3N/4} = \probawithstartingpointmodified{x}{\absolutevalue{R(X_{4m})} + \absolutevalue{B(X_{4m})} \geq N/4} \leq \epsilon/32.\]
    We conclude that 
    \begin{align*}
        \probawithstartingpointmodified{x}{L_{4m} \geq N/2} &\geq \probawithstartingpointmodified{x}{\absolutevalue{G(X_{2m})} < N/4, \absolutevalue{G(X_{4m})} > 3N/4} \\
        &\geq 1 - \probawithstartingpointmodified{x}{\absolutevalue{G(X_{2m})} \geq N/4} - \probawithstartingpointmodified{x}{\absolutevalue{R(X_{4m})} \leq 3N/4}\\
        &\geq 1 - \epsilon/16.    
    \end{align*}
    This implies that, in the original dynamics, we also have
    \[ \probawithstartingpoint{x}{L_{4m} \geq N/2} \geq 1 - \epsilon/16.\]
    Therefore, by an argument of union bound and the Markov property, 
    \[\probawithstartingpoint{x}{L_{16m} \geq 2N} \geq 1 - 4\times \epsilon/16 = 1- \epsilon/4.\]
    We choose $t_2 = 16m$ to conclude the proof.
\end{proof}
\begin{remark}
    The number $t_2$ above is not meant to be optimal. In fact, if we are interested only in the case of large $N$, we can choose $m = \dfrac{1+ o(1)}{c} \log 4$. 
\end{remark}
The rest of this section is devoted to proving the anticoncentration inequality in Lemma \ref{lemma_anticoncentration_inequality}. We introduce some notations that we will use in Lemma \ref{lemma_anticoncentration_inequality} and Proposition \ref{proposition_crossing_inequality}. Let $(T_i)_{i \geq 0}$ be defined as in Observation \ref{observation_conditional_evolution}. We fix $t \in \Rbb_+$ and denote by $T_{a_1},\dots, T_{a_r}$ the times at which a green individual is recolored blue at site $1$, and $T_{b_1}, \dots,T_{b_s}$ the times at which a red or blue individual is recolored green at site $N$, during the time interval $[0,t]$. All these times are $\greentrajectory$-measurable. By an abuse of notation, we relabel those individuals by $\tilde{a}_1, \dots, \tilde{a}_r, \tilde{b}_1, \dots, \tilde{b}_s$. We denote the random walks killed when being recolored at site $N$ of the individuals $\tilde{a}_1, \dots, \tilde{a}_r$ by $A_1,\dots, A_r$, respectively. If the individual labelled $i$ is blue or red at time $0$, we denote by $\sigmatilde(i,\cdot)$ the walk of the individual $i$ killed when being recolored at site $N$ as well. Observe that the walk $\sigmatilde(i, \cdot)$ either survives up to time $t$ or is killed at some time $T_{b_l}$. We will use the notation $\{\sigmatilde(i, t) = j\}$ to indicate the event that the walk survives up to time $t$ and ends up at site $j$, and $\{\sigmatilde(i, t) = \tilde{b}_l\}$ to say that the walk is killed at time $T_{b_l}$. Similar notations are used for the walk $A_k, \, 1\leq k \leq r$\\
We will need the following proposition.
\begin{proposition}[Crossing inequality]\label{proposition_crossing_inequality}
    For any $x \in \setcolorinterchange$, for any individual labelled $i$ that is either blue or red in $x$, for any site $j$, for any $1 \leq k \leq r, \, 1 \leq l \leq s$,
    \begin{equation}\label{crossing_inequality}
        \probawithstartingpoint{x}{\sigmatilde(i,t) = j, A_k(t) = \tilde{b}_l \big| \greentrajectory} \leq\probawithstartingpoint{x}{\sigmatilde(i,t) = \tilde{b}_l, A_k(t) = j \big| \greentrajectory}.
    \end{equation}
\end{proposition}
We show how we can prove Lemma \ref{lemma_anticoncentration_inequality} using Proposition \ref{proposition_crossing_inequality}.
\begin{proof}[Proof of Lemma \ref{lemma_anticoncentration_inequality}]
Let $i, j$ be as in the statement of Lemma \ref{lemma_anticoncentration_inequality}. Summing the inequality in Proposition \ref{proposition_crossing_inequality} over $l$, we get 
    \begin{equation*}
        \probawithstartingpoint{x}{\sigmatilde(i,t) = j, \text{$A_k$ is killed by time $t$} \big| \greentrajectory} \leq \probawithstartingpoint{x}{\text{$\sigmatilde(i,\cdot)$ is killed by time $t$}, A_k(t) = j\big| \greentrajectory}.
    \end{equation*}
    Now summing over $k$, we get
    \begin{equation}\label{label1}
    \begin{split}
        &\esperancewithstartingpoint{x}{\indicator{\sigmatilde(i,t) = j} \times \#\{\text{walks among $A_1,\dots,A_r$ that are killed by time $t$}\} \big| \greentrajectory}\\
        &\leq \sum_{k=1}^s \probawithstartingpoint{x}{\text{$\sigmatilde(i,\cdot)$ is killed by time $t$}, A_k(t) = j \big| \greentrajectory}.
    \end{split}    
    \end{equation}
     On the event $\{L_t \geq 2N\}$, there are at least $2N$ times at which a blue or red individual is recolored green by time $t$. Since originally there are at most $N$ blue and red individuals, then there are always at least $N$ times a green individual is recolored blue at site $1$ and then recolored green again by time $t$. In other words, at least $N$ walks are killed at site $N$ by time $t$ among $A_1,\dots, A_r$. Hence the left-hand side of equation \eqref{label1} is at least $N\times \proba{\sigmatilde(i,t) = j\big| \greentrajectory}$. On the other hand, we can realize the trajectories of $\sigmatilde(i,\cdot)$ and $A_1, A_2,\dots, A_r$ altogether by revealing $\Xi^{BR}$. Subsequently, the events on the right-hand side of \eqref{label1} are pairwise disjoint since $A_k$ and $A_l$ cannot both occupy site $j$ at time $t$, for any $k \neq l$. This implies that the sum on the right-hand side of \eqref{label1} is smaller than $1$. So overall, on the event $\{L_t \geq 2N\}$, 
     \begin{equation}\label{anticoncentration}
        N\times \probawithstartingpoint{x}{\sigmatilde(i,t) = j \big| \greentrajectory} \leq 1.
     \end{equation}
    Furthermore, 
    \begin{equation}\label{anticoncentration_equality}
        \probawithstartingpoint{x}{\sigmatilde(i,t) = j \big| \greentrajectory} = \probawithstartingpoint{x}{\sigma_t(i) = j, \sigma_s(i) \notin G(X_s), \forall0 \leq s \leq t\big| \greentrajectory}.
    \end{equation}
    Combining \eqref{anticoncentration} and \eqref{anticoncentration_equality}, we get what we want.
\end{proof}
We now prove Proposition \ref{proposition_crossing_inequality}.
\begin{proof}[Proof of Proposition \ref{proposition_crossing_inequality}]
    The idea is to use the monotonicity of the exclusion process. By Observation \ref{observation_conditional_evolution}, we see that $\sigmatilde(i, \cdot)$ is a simple random walk conditionally on the environment created by the green region. If that walk is still alive at the time when $\tilde{a}_k$ is born, then it must be on the right of $\tilde{a}_k$ at that time (since $\tilde{a}_k$ is born at site $1$). Those two individuals then evolve as a simple exclusion process with two individuals conditionally on the environment as we realize $\Xi^{BR}$, still by Observation \ref{observation_conditional_evolution}. We propose another graphical construction of those two walks as follows. At any time $t$, if there are two walks alive, say at two sites $u, v$ with $u < v$, we refer to the individual at site $u$ as the individual on the left and to the individual at site $v$ as the individual on the right. In case there is only one walk alive, we refer to it as the only individual. Let $\Xi_R, \Xi_L, \Xi_E$ be $3$ independent Poisson processes with $\Xi_R$ and $\Xi_L$ of intensity $N^2\drm t \otimes \textrm{Card}$ on $\Rbb_+ \times \{L, R\}$, and $\Xi_E$ of intensity $2N^2 \drm t$. R, L, and E stand for right, left, and exchange. The rule of evolution is as follows. 
    \begin{itemize}
        \item For each point $(s, w)$ of $\Xi_R$, the individual on the right at time $s-$ (or the only individual if there is only one walk alive) attempts to make a jump at time $s$ to the left if $w = L$ or to the right if $w = R$. It succeeds except when trying to jump on a site occupied by a green individual or the individual on the left, in which case the jump is cancelled. 
        \item For each point $(s, w)$ of $\Xi_L$, the individual on the left at time $s-$ (if there exists) attempts to make a jump at time $s$ to the left if $w = L$ or to the right if $w = R$. It succeeds except when trying to jump on a site occupied by a green individual or by the individual on the right, in which case the jump is cancelled.
        \item For each point $s$ of $\Xi_E$, we exchange the positions of the two individuals with probability $1/2$ if they are adjacent at time $s$. 
        \item For each time $T_i$, the two walks make the jump forced by the environment, as explained in Observation \ref{observation_conditional_evolution}.
    \end{itemize}
    This construction gives us the same distribution of the two walks as the one given by $\Xi^{BR}$, as the rates of transition given by the two constructions are always the same. In short, $\Xi_R$ and $\Xi_L$ are used to generate the walks on the right and the left, respectively. When we realize $\Xi_R, \Xi_L$, we observe two random walks that are not allowed to jump on top of each other. $\Xi_E$ is then used to make precise where those two walks exchange their positions. 
    With this construction, the set $\{\sigmatilde(i, s), A_k(s)\}$ is entirely determined by $\Xi_R, \Xi_L$, for any $0 \leq s \leq t$, since $\Xi_E$ has no effect on that set. Now conditionally on the $(\Xi_R, \Xi_L)$-measurable event $\{\sigmatilde(i, t), A_k(t)\} = \{j, \tilde{b}_l\}$ (which means that one walk reaches $j$ at time $t$ and the other walk was killed at $b_l$), we realize $\Xi_E$ to obtain a number $m$ of exchange edges between the two walks. We see that if $m=0$, then we necessarily have $\sigmatilde(i,t) = \tilde{b}_l, A_k(t) = i$ due to monotonicity, which means
   \begin{equation}\label{crossing_when_m_is_zero}
        0 = \probawithstartingpoint{x}{\sigmatilde(i,t) = j, A_k(\timeincrossinglemma) = \tilde{b}_l, m = 0 \big| \greentrajectory} \leq \probawithstartingpoint{x}{\sigmatilde(i,t) = \tilde{b}_l, A_k(\timeincrossinglemma) = j, m = 0 \big| \greentrajectory}.
    \end{equation}
   If $m > 0$, we cannot distinguish the two individuals anymore since the probability that the sum of $m$ independent Bernoulli variables of parameter $1/2$ is even (or odd) is $1/2$, which means
    \begin{equation} \label{crossing_when_m_is_positive}
        \probawithstartingpoint{x}{\sigmatilde(i,t) = j, A_k(\timeincrossinglemma) = \tilde{b}_l, m > 0 \big| \greentrajectory} = \probawithstartingpoint{x}{\sigmatilde(i,t) = \tilde{b}_l, A_k(\timeincrossinglemma) = j, m > 0 \big| \greentrajectory}, 
        %= 1/2 \proba{\{X_t(\timeincrossinglemma), A_k(\timeincrossinglemma)\} = \{j, b_l\}, m > 0} 
    \end{equation}
    since both are equal to $\dfrac{1}{2} \probawithstartingpoint{x}{\{X_t(\timeincrossinglemma), A_k(\timeincrossinglemma)\} = \{j, \tilde{b_l}\}, m > 0 \big| \greentrajectory}$. Summing \eqref{crossing_when_m_is_zero} and \eqref{crossing_when_m_is_positive}, we get what we want.
\end{proof}

\section{The lower bound}
We finish the proof of Theorem \ref{maintheoremtwo} by proving the lower bound on $\tmix$ in Theorem \ref{maintheoremtwo}.
\begin{proof}[Proof of the lower bound]
    We follow Wilson's classical method in \cite{Wilson2004}. We consider the process with the initial condition $\eta_0 = \mathds{1}$. Recall that the weight of a configuration $\eta$ is 
    $$ S(\eta) := \sum_{i = 1}^N \eta(i).$$
    $S(\eta)$ will be our distinguishing statistic. Recall that, for any $t \geq 0$,
    \[\eta_t \overset{d}{=} \pushforwardfunction(\eta, X_t, \xi^B, \xi^G),\]
    with $X_t, \xi^B, \xi^G$ defined as in Lemma \ref{natural_coupling}.
    We write $S_t$ for $S\left(\pushforwardfunction(\eta, X_t, \xi^B, \xi^G)\right)$, and $S_\infty$ for $S(\eta_\infty)$, where $\eta_\infty \sim \pi$. By symmetry, we see that 
    \begin{equation*}
        \esperance{\absolutevalue{B(X_t)}} = \esperance{\absolutevalue{G(X_t)}} = \dfrac{N -\esperance{\absolutevalue{R(X_t)}}}{2}.
    \end{equation*}
    Hence 
    \begin{align*}
        \esperance{S_t} &= \esperance{\esperance{S_t\big|X_t}}\\
        &= \esperance{\absolutevalue{R(X_t)} + p\absolutevalue{B(X_t)} + q\absolutevalue{G(X_t)}}\\
        &= \esperance{\absolutevalue{R(X_t)}} + \dfrac{N - \esperance{\absolutevalue{R(X_t)}}}{2}(p+q),
    \end{align*}
    Letting $t \to \infty$, and observing that $\esperance{\absolutevalue{R(X_t)}} \xrightarrow{ t \to \infty} 0$ by Lemma \ref{exponential_decay_w0}, we get
    \begin{align*}
        \esperance{S_\infty} &= N \times \dfrac{p+q}{2},\\
        \esperance{S_t} - \esperance{S_\infty} &= \esperance{\absolutevalue{R(X_t)}} \left(1 - \dfrac{p+q}{2}\right).
    \end{align*}
    By Proposition \ref{strongly_rayleigh_preserved_by_exclusion}, $\eta_t$ is ND. Hence by inequality \eqref{sum_of_strongly_Rayleigh_vector},
    \[\variance{S_t} \leq \esperance{S_t},\]
    and by letting $t \to \infty$, 
    \[\variance{S_\infty} \leq \esperance{S_\infty}.\]
    Hence, under assumption \eqref{assumption}, 
    \[\maxset{\variance{S_t}, \variance{S_\infty}} \leq \dfrac{N(p+q)}{2} + \esperance{\absolutevalue{R(X_t)}} \left(1 - \dfrac{p+q}{2}\right) \leq Np + \esperance{\absolutevalue{R(X_t)}}.\]
    So by Proposition $7.9$ in \cite{Levin2017}, 
    \[\norm{P^t_{\eta_0} - \pi}_{TV} \geq 1 - 8 \dfrac{\maxset{\variance{S_t}, \variance{S_\infty}}}{(\esperance{S_t} - \esperance{S_\infty})^2}  \geq 1- 8 \dfrac{Np + \esperance{\absolutevalue{R(X_t)}}}{\left(1 - \dfrac{p+q}{2}\right)^2 \esperance{\absolutevalue{R(X_t)}}^2}.\]
    Let $t = t^* - 2m$, for some positive integer $m < t^*/2$ that we will choose later, and let $c$ be the constant in Lemma \ref{one_walk_marginal}. By Lemma \ref{exponential_decay_w0} and the Markov property, $\esperance{\absolutevalue{R(X_t)}} \geq \esperance{\absolutevalue{R(X_{t^*})}} e^{cm} =  (\sqrt{Np} \vee 1)e^{cm}$, hence the last term is bigger than 
    \[1 - \dfrac{8}{\left(1 - \dfrac{p+q}{2}\right)^2 e^{2cm}} - \dfrac{8}{\left(1 - \dfrac{p+q}{2}\right)^2 e^{cm}}.\]
    Note that with our assumption, $1 > 1 -\dfrac{p+q}{2} > \dfrac{1}{2}$. So we conclude that the term above is greater than $1 - \dfrac{32}{e^{2cm}} - \dfrac{32}{e^{cm}}$, which is greater than $\epsilon$ for $m = \left\lceil\dfrac{1}{c} \log\left(\dfrac{64}{1-\epsilon}\right)\right\rceil$. This implies that $\tmix\geq t^* - 2m$, if $m <t^*/2$. Besides, this inequality is trivial when $m \geq t^*/2$, which finishes our proof.
\end{proof}
\begin{remark}
    This method can also give a lower bound on the mixing time from any initial configuration $\eta$. In fact, we can prove that, for some constant $C$, 
    \[\tmixetaepsilon \geq t^*(\eta) - C \left(1 + \log\left(\dfrac{1}{1-\epsilon}\right)\right),\]
    where 
    \[t^*(\eta) = \minset{t \geq 0: \absolutevalue{\Ebb_{\eta}^0\left[S_t\right] - \esperance{S_{\infty}}} \leq \sqrt{Np} \vee 1}. \]
\end{remark}
\section{The computation of $t^*$}
For the sake of completeness, we include here the proof of equality \eqref{estimation_of_tmix}, which is a particular case of the computations presented in Appendix A of \cite{Goncalves2021}.
\begin{proof}[Proof of \eqref{estimation_of_tmix}]
    We consider the model where the reservoirs are of density $p = q = 0$. We define $u_t: [N] \to [0,1]$ by $u_t(x) = \esperancewithstartingpoint{\mathds{1}}{\eta_t(x)}$. Then by Dynkin's formula, $\{u_t; \, t\geq 0\}$ is the unique solution of the system of equations 
    \begin{equation*}
        \left\lbrace
            \begin{array}{@{}r@{\;}l@{}l@{}}
                \ddt u_t(x) &= \Delta u_t(x) &\text{\quad for $t\geq 0$ and $x \in [N]$} , \\
                 u_0(x) &= 1 &\text{\quad for $x \in [N]$} ,
            \end{array}
        \right.
    \end{equation*}
    where $\Delta$ is the discrete Laplacian defined by, for any function $f: [N] \to \Rbb$, 
    \[\Delta f = N^2 (f(x+1) + f(x-1) - 2f(x)), \, \forall x \in [N],\]
    with the convention that $f(0) = f(N+1) = 0$. The eigenfunctions of $\Delta$ are given by 
    \[\phi_l(x) = \sqrt{2} \sin \left(\dfrac{\pi l x}{N+1}\right), \quad 1\leq l \leq N,\]
    with the corresponding eigenvalues $-\lambda_l$ given by 
    \[\lambda_l = 2N^2\left(1 - \cos \left(\dfrac{\pi l }{N+1}\right)\right), \quad 1\leq l \leq N.\] 
    Besides, $(\phi_l)_{ 1\leq l \leq N}$ is an orthonormal basis for the scalar product given by 
    \[\crochet{f}{g} = \dfrac{1}{N+1} \sum_{x = 1}^N f(x)g(x),\]
    for any $f, g:\, [N] \to \Rbb$. Let $t$ be an arbitrary positive number. We see that
    \begin{align}\label{numberofparticlescalarproduct}
        \esperancewithstartingpoint{\mathds{1}}{\absolutevalue{R(X_t)}} &= \sum_{x = 1}^N u_t(x) = (N+1) \crochet{u_t}{u_0}.
    \end{align}
    Let 
    \[u_0 = \sum_{l = 1}^N c_l \phi_l\]
    be the decomposition of $u_0$ in the basis $(\phi_l)_{1 \leq l \leq N}$. Then 
    \begin{align*}
        u_t = \sum_{l = 1}^N c_l e^{-\lambda_l t} \phi_l.
    \end{align*}
    Hence 
    \begin{align}
        \crochet{u_t}{u_0} = \sum_{l = 1}^N c_l^2 e^{-\lambda_l t}.
    \end{align}
    Note that $\lambda_1 = \maxset{\lambda_1, \dots, \lambda_N}$. Hence 
    \begin{align*}
        c_1^2 e^{-\lambda_1 t} \leq \crochet{u_t}{u_0} \leq (\sum_{l = 1}^N c_l^2) e^{-\lambda_1 t} = \crochet{u_0}{u_0} e^{-\lambda_1 t} = \dfrac{N}{N+1} e^{-\lambda_1 t}.
    \end{align*}
    Plug in \eqref{numberofparticlescalarproduct}, we see that
    \begin{align*}
        (N+1) c_1^2 e^{-\lambda_1 t} \leq \esperancewithstartingpoint{\mathds{1}}{\absolutevalue{R(X_t)}} \leq Ne^{- \lambda_1 t}, 
    \end{align*}
    which means 
    \begin{align}\label{inequalityt*}
        \dfrac{1}{\lambda_1} \left(2\log (c_1) + \log \left(\dfrac{N+1}{\esperancewithstartingpoint{\mathds{1}}{\absolutevalue{R(X_t)}}}\right)\right)\leq t \leq \dfrac{1}{\lambda_1} \log \left(\dfrac{N}{\esperancewithstartingpoint{\mathds{1}}{\absolutevalue{R(X_t)}}}\right).
    \end{align}
    We finish the proof by estimating $\lambda_1$ and $c_1$.
    Note that by Taylor's expansion of function cosine around $0$, 
    \begin{equation}\label{lambdaestimation}
        \lambda_1 = 2N^2 \left(\dfrac{1}{2} \dfrac{\pi^2}{(N+1)^2} + \Ocal{\dfrac{\pi^4}{(N+1)^4}}\right) = \pi^2 + \Ocal{1/N}.    
    \end{equation}
    Besides, 
    \begin{align*}
        c_1 &= \crochet{\mathds{1}}{\phi_1}\\
        &= \dfrac{1}{N+1} \sum_{x = 1}^N \phi_1(x) \\
        &= \dfrac{\sqrt{2}}{N+1} \sum_{x = 1}^N \sin \dfrac{\pi x}{N+1}.
    \end{align*}
    By some classical trigonometric computations, 
    \begin{align}\label{c_one_estimation}
        c_1 = &\dfrac{\sqrt{2}}{N+1} \times \dfrac{\cos \dfrac{\pi}{2(N+1)}}{\sin \dfrac{\pi}{2(N+1)}} = \dfrac{2\sqrt{2}}{\pi} \left(1 + \Ocal{1/N^2}\right),
    \end{align}
    where we have used the Taylor expansion of the sinus and cosinus functions around $0$. Note that at $t^*$, $\esperancewithstartingpoint{\mathds{1}}{\absolutevalue{R(X_{t^*})}} = \sqrt{Np} \vee 1$, by \eqref{t*andredregion}. Then the three equations \eqref{inequalityt*}, \eqref{lambdaestimation}, \eqref{c_one_estimation} together imply
    \begin{equation*}
         t^* = \dfrac{1}{\pi^2} \log \left(\dfrac{N}{\sqrt{Np} \vee 1}\right) \pm \Ocal{1},
    \end{equation*}
    which is precisely what we want.
\end{proof}    

\bibliographystyle{abbrv}
\bibliography{ref}

\end{document}